\documentclass[12pt]{amsart}

\newcommand{\newsection}[1]{\setcounter{equation}{0} \section{#1}}
\newcommand{\bea}{\begin{eqnarray}}
\newcommand{\eea}{\end{eqnarray}}

\newcommand{\clb}{\mathcal{B}}
\newcommand{\clc}{\mathcal{C}}
\newcommand{\cld}{\mathcal{D}}
\newcommand{\cle}{\mathcal{E}}
\newcommand{\clf}{\mathcal{F}}

\newcommand{\clh}{\mathcal{H}}
\newcommand{\clk}{\mathcal{K}}
\newcommand{\cll}{\mathcal{L}}
\newcommand{\clm}{\mathcal{M}}
\newcommand{\cln}{\mathcal{N}}
\newcommand{\clo}{\mathcal{O}}

\newcommand{\clx}{\mathcal{X}}
\newcommand{\cly}{\mathcal{Y}}

\newcommand{\raro}{\rightarrow}

\def \qed {\hfill \vrule height6pt width 6pt depth 0pt}
\def\textmatrix#1&#2\\#3&#4\\{\bigl({#1 \atop #3}\ {#2 \atop #4}\bigr)}
\def\dispmatrix#1&#2\\#3&#4\\{\left({#1 \atop #3}\ {#2 \atop #4}\right)}
\newcommand{\be}{\begin{equation}}
\newcommand{\ee}{\end{equation}}
\newcommand{\ben}{\begin{eqnarray*}}
\newcommand{\een}{\end{eqnarray*}}

\newcommand{\NI}{\noindent}

\newcommand{\bi}{\begin{itemize}}
\newcommand{\ei}{\end{itemize}}

\def\5{{5\superprime}}

\newtheorem{Theorem}{\sc Theorem}[section]
\newtheorem{Lemma}[Theorem]{\sc Lemma}
\newtheorem{Proposition}[Theorem]{\sc Proposition}
\newtheorem{Corollary}[Theorem]{\sc Corollary}
\newtheorem{Definition}[Theorem]{\sc Definition}
\newtheorem{Example}[Theorem]{\sc Example}
\newtheorem{Remark}[Theorem]{\sc Remark}
\newtheorem{Note}[Theorem]{\sc Note}
\newtheorem{Question}[Theorem]{\sc Question}
\newtheorem{ass}[Theorem]{\sc Assumption}
\newcommand{\bt}{\begin{Theorem}}
\def\beginlem{\begin{Lemma}}
\def\beginprop{\begin{Proposition}}
\def\begincor{\begin{Corollary}}
\def\begindef{\begin{Definition}}
\def\beginexamp{\begin{Example}}
\def\beginrem{\begin{Remark}}
\def\beginq{\begin{Question}}
\def\beginass{\begin{ass}}
\def\beginnote{\begin{Note}}
\newcommand{\et}{\end{Theorem}}
\def\endlem{\end{Lemma}}
\def\endprop{\end{Proposition}}
\def\endcor{\end{Corollary}}
\def\enddef{\end{Definition}}
\def\endexamp{\end{Example}}
\def\endrem{\end{Remark}}
\def\endq{\end{Question}}
\def\endass{\end{ass}}
\def\endnote{\end{Note}}

\begin{document}

\title[On c.n.c. commuting contractive tuples]{On c.n.c. commuting contractive tuples}

\author[Bhattacharyya]{T. Bhattacharyya}

\address{Department of Mathematics, Indian Institute of Science,
Bangalore 560012, India}

\email{tirtha@math.iisc.ernet.in}

\author[Eschmeier]{J. Eschmeier}

\address{Fachbereich Mathematik, Universit\"{a}t des Saarlandes, 66123
Saarbr\"{u}cken, Germany}

\email{eschmei@math.uni-sb.de}

\author[Sarkar]{J. Sarkar}

\address{Department of Mathematics, Indian Institute of Science,
Bangalore 560012, India}

\email{jaydeb@math.iisc.ernet.in}

\begin{abstract}
The characteristic function has been an important tool for studying
completely non unitary contractions on Hilbert spaces. In this note,
we consider completely non-coisometric contractive tuples of
commuting operators on a Hilbert space $\clh$. We show that the
characteristic function, which is now an operator valued analytic
function on the open Euclidean unit ball in $\mathbb{C}^n$, is a
complete unitary invariant for such a tuple. We prove that the
characteristic function satisfies a natural transformation law under
biholomorphic mappings of the unit ball. We also characterize all
operator-valued analytic functions which arise as characteristic
functions of pure commuting contractive tuples.
\end{abstract}

\maketitle

\newsection{Introduction}

The characteristic function for a single contraction on a Hilbert
space was defined by Sz.-Nagy and Foias in \cite{S-N}. Since then it
has drawn a lot of attention and several interesting results are
known about it.

A tuple $T = (T_1, \ldots ,T_n)$ of bounded operators on a Hilbert
space $\clh$ is called contractive if $\|T_1 h_1 + \cdots + T_n h_n
\|^2 \le \|h_1 \|^2 + \cdots + \|h_n\|^2$ for all $h_1, \ldots ,
h_n$ in $\mathcal{H}$ or equivalently $\sum_{i=1}^n T_iT_i^* \le
1_{\mathcal{H}}$. The positive operator $ (1_{\mathcal{H}} -
\sum_{i=1}^n T_iT_i^*)^{1/2}$ and the closure of its range will be
called the {\em defect operator} $D_{T^*}$ and the {\em defect
space} $\cld_{T^*}$ of $T$.

We shall also denote by $T$ the row operator from $\clh^n$ to $\clh$
which maps $(h_1, \ldots , h_n)$ to $ T_1h_1 + \cdots + T_nh_n$. The
adjoint $T^* : \clh \raro \clh^n$ maps $h$ to the column vector
$(T_1^*h, \ldots , T_n^*h)$ and, in fact, $T$ is a contractive tuple
if and only if the operator $T$ is a contraction. Thus for a
contractive tuple $T$ one can also consider the defect operator $D_T
= (1_{\clh^n} - T^*T)^{1/2}$ $= (( \delta_{ij} 1_\clh - T_i^* T_j
))^{1/2}$ in $\clb(\clh^n)$ and the associated defect space $\cld_T
= \overline{\rm Ran} D_T \subset \clh^n$.

We use the notation $ \mathbb{B}_n$ for the open Euclidean unit ball
in $ \mathbb{C}^n$. The prototypical example, which has been used by
Arveson \cite{sub3}, M\"uller and Vasilescu \cite{MV} in the
construction of appropriate models, is the shift on $H^2_n$ defined
as follows. Given a complex Hilbert space $\cle$, let $\clo (
\mathbb{B}_n , \cle)$ be the class of all $\cle$-valued analytic
functions on $\mathbb{B}_n$. For any multi-index $k = (k_1, \ldots
,k_n) \in \mathbb{N}^n$, we write $|k| = k_1 + \cdots +k_n$. Then
consider the Hilbert space \bea H^2_n(\cle) = \{ f = \sum_{k \in
\mathbb{N}^n} a_k z^k \in \clo( \mathbb{B} , \cle) : \; a_k \in \cle
\; \mbox{with} \; \| f \|^2 = \sum_{k \in \mathbb{N}^n} \frac{\| a_k
\|^2}{ \gamma_k} < \infty \}, \label{HEdef} \eea where $\gamma_k =
|k | ! / k !$. One can show that $H^2_n(\cle)$ is the $\cle$-valued
functional Hilbert space given by the reproducing kernel $ ( 1 -
\langle z , w \rangle )^{-1} 1_\cle$. When $\cle = \mathbb{C}$, we
use the abbreviation $H^2_n$. For $n=1$, this space is the usual
Hardy space on the unit disk. The space $H^2_n(\cle)$ is
isometrically isomorphic to the Hilbertian tensor product $H^2_n
\otimes \cle$ in a canonical way. Given complex Hilbert spaces
$\cle$ and $\cle_*$, the multiplier space $ M(\cle , \cle_*)$
consists of all $\varphi \in \clo( \mathbb{B}_n , \clb(\cle ,
\cle_*))$ such that $\varphi H^2_n(\cle) \subset H^2_n(\cle_*) .$ By
the closed graph theorem, for each function $\varphi \in M(\cle ,
\cle_*)$,the induced multiplication operator $M_\varphi : H(\cle)
\raro H(\cle_*),$ $\, f \mapsto \varphi f$, is continuous. The {\em
standard shift} on $H^2_n(\cle)$ is the tuple $M^\cle_z =
(M^\cle_{z_1}, \ldots ,M^\cle_{z_n})$ consisting of the
multiplication operators $M^\cle_{z_i} : H^2_n(\cle) \rightarrow
H^2_n(\cle)$ with the coordinate functions $z_i$. When $\cle =
\mathbb{C}$, we shall write $M_z$ for $M_z^\cle$. Arveson studied
both the space $H^2_n$ and the standard shift $M_z$ in great detail
in \cite{sub3}. It can be seen without much difficulty that $M_z$ is
a commuting contractive tuple. In fact, $D_{M_z^*}$ is the
one-dimensional projection onto the space of constant functions. The
space $H^2_n$ was first used by Drury \cite{Drury} who generalized
von Neumann's inequality to operator tuples.

With a commuting contractive tuple $T$ on a Hilbert space
$\mathcal{H}$, one associates a completely positive map $ P_T
:\mathcal{B}(\mathcal{H})
 \raro \mathcal{B}(\mathcal{H}) $ defined by
$P_T(X)= \sum_{i=1}^n {T_i}XT_i^*$. We denote by $A_T \in \clb
(\clh)$ the strong limit of the decreasing sequence of positive
operators $I \geq P_T(I) \geq P_T^2(I) \geq... \geq 0$. The tuple
$T$ is called pure if $A_T = 0$. For $n=1$, this corresponds to the
$C_{\cdot 0}$ case in the Sz.-Nagy and Foias classification
\cite{S-N} of contractions. The standard shift $M_z$ on $H^2_n$ is
pure.

An operator valued bounded function on $\mathbb{B}_n$ is a triple
$\{ \cle , \cle_* , \varphi\}$, where $\cle$ and $\cle_*$ are
Hilbert spaces and $\varphi$ is a $\clb( \cle , \cle_*)-$valued
bounded function on $\mathbb{B}_n$. If $\| \varphi (z) \| \le 1$,
then the function is called contractive. Two operator valued bounded
functions $\{ \cle , \cle_* , \varphi\}$ and $\{ \clf , \clf_* ,
\psi \}$ are said to coincide if there exist unitary operators $\tau
: \cle \raro \clf$ and $ \tau_* : \cle_{*} \raro \clf_{*}$ such that
the following diagram

 \setlength{\unitlength}{3mm}
 \begin{center}
 \begin{picture}(20,14)(0,0)
 \put(2,3){$ \clf$} \put(10,3){$ \clf_{*}$}
 \put(5.6,2.2){$ \psi(z)$}
 \put(1.5,6.5){$ \tau$} \put(11,6.5){$ \tau_*$}
 \put(2,10){$ \cle$} \put(10,10){$ \cle_{*}$}
 \put(5.6,11){$ \varphi(z)$}
 \put(3.5,3.5){ \vector(1,0){6}} \put(3.5,10.5){ \vector(1,0){6}}
 \put(2.4,9.2){ \vector(0,-1){5}} \put(10,9.2){ \vector(0,-1){5}}

 \end{picture}
 \end{center}

commutes for all $z$ in $ \mathbb{B}_n$.

The characteristic function for a commuting contractive tuple is
defined as the operator valued contractive function $ \{ \cld_T ,
\cld_{T^*} , \theta_T \}$, where \bea \theta_T (z) = - T + D_{T^*} (
1_\clh - Z T^* )^{-1} Z D_T, \;\; z \in \mathbb{B}_n.
\label{chfndef} \eea Here $Z = (z_1 I_{\clh}, \ldots ,z_n I_{\clh})$
denotes the row multiplication induced by $z \in \mathbb{B}_n.$ The
characteristic function was defined in \cite{BES}, and it was proved
to be a complete unitary invariant in the case of pure tuples $T$.
Theorem 4.4 in  \cite{BES} states that if $T$ and $R$ are two pure
commuting contractive tuples on Hilbert spaces $\clh$ and $\clk$,
respectively, then $T$ and $R$ are unitarily equivalent (that is,
$T_i = U R_i U^*$ for all $i=1, \ldots , n$ and a suitable unitary
operator $U : \clk \raro \clh$) if and only if $ \{ \cld_T ,
\cld_{T^*} , \theta_T \}$ and $ \{ \cld_R , \cld_{R^*} , \theta_R
\}$ coincide. We would like to point out here that recent works by
Gelu Popescu also shows the same result in a more general setting
(see \cite{Pop05a} and \cite{Pop05b}). Popescu has defined the
characteristic function for a contractive tuple in \cite{Pop89}. He
proved that it is a complete unitary invariant for  completely
non-coisometric tuples. Since he considered the non- commutative
case, his characteristic function was defined as a multi analytic
operator. In \cite{Pop05a} he has shown that one can associate a
constrained characteristic function to a constrained contractive
tuple (a commuting contractive tuple is a particular example). Thus
he also obtained the result mentioned above by compressing the
multi-analytic operator to the symmetric Fock space.

This note serves the purpose of proving three basic results about
the characteristic function. We show that the characteristic
function is a complete unitary invariant for completely
non-coisometric commuting contractive tuples. In the process, we
construct a functional model for a completely non-coisometric tuple.
This is the content of Section 3. We start in Section 2 by showing
that the characteristic function obeys a certain natural
transformation rule with respect to automorphisms of the Euclidean
unit ball. The automorphisms of the disk have played an important
role in the model and dilation theory of single contractions. Hence
it is naturally desirable to obtain a multivariable analogue.
Finally, in Section 4, we characterize the subspaces of $H^2_n
\otimes \cle$ reducing for the canonical shift and use this result
to describe all operator-valued analytic functions which arise as
characteristic functions of pure tuples. We also prove a version of
the classical Beurling-Lax-Halmos theorem in that section using the
characteristic function.

\vspace{1mm}

{\small After completion of this note, we came to know about a
recent preprint of Benhida and Timotin \cite{BTim} which also
studies the connection between commuting contractive tuples and
automorphisms of the unit ball.}

\newsection{Transformation Rule}

For a fixed $a \neq 0$ in $\mathbb{B}_n$, let $P_a$ denote the
orthogonal projection of $\mathbb{C}^n$ onto the one dimensional
subspace $[ a ]$ generated by $a$, that is,  $P_a z =  (\langle z, a
\rangle / \langle a, a \rangle)  a$. Let $Q_a$ be the orthogonal
projection $I - P_a $. Define $s_a = ( 1 - |a|^2)^{1/2}$. Then  $$
\varphi_a (z) = \frac{a - P_a z - s_a Q_a z}{1 - \langle z , a
\rangle}$$ is a biholomorphic automorphism of the unit ball (see
\cite{Rudin}). Given a commuting contractive tuple of operators $T$,
it is easy to see that the Taylor spectrum of $T$ is contained in
the closure of the unit ball. Since $\varphi_a$ is actually analytic
in an open set containing the closed unit ball, we can consider the
associated operator tuple
$$T_a = (1 - T A^*)^{-1} ( A - P_a T - s_a Q_a T),$$
where $A = (a_1 I_\clh , \ldots , a_n I_\clh)$, $ P_a T = \frac{1}{
|a|^2} T A^* A$ and $Q_a T = T - P_a T$.

\begin{Lemma}

Let $T$ be a commuting contractive $n$-tuple of operators on a
Hilbert space $\clh$. Then for any $a \in \mathbb{B}_n$, we have the
identiy
$$
I_\clh - T_a T_a ^* = ( 1 - |a|^2) (I_\clh - T A^*)^{-1} ( I_\clh -
T T^*) (I_\clh - A T^*)^{-1}.
$$

\end{Lemma}

\vspace{0.3cm}

\NI \textsf{Proof.} Using the equality  $ P_a T (Q_a T)^*  = 0$ and
the definition of $T_a$, we obtain that \ben T_a T_a ^* & = & (I - T
A^*)^{-1} ( A - P_a T - s_a Q_a T) ( A^* - (P_a T)^* - s_a (Q_a T)^*
) (I - A
T^*)^{-1} \\
& = & (I - T A^*)^{-1} [ |a|^2 I - A T^* - T A^* + \frac{1}{|a|^2} T
A^* A T^* + s_a ^2 ( T - \frac{1}{|a|^2} T A^* A) \\
& & ( T^* - \frac{1}{|a|^2} A^* A T^*)] (I - A T^*)^{-1} \\
& = & (I - T A^*)^{-1} [ |a|^2 I - A T^* - T A^* + T A^* A T^* + (1
- |a|^2) T T^* ] (I - A T^*)^{-1} \\ & = & (I - T A^*)^{-1} [ (I - T
A^*) (I - A T^*) - (1 - |a|^2) + (1 - |a|^2) T T^* ] (I - A
T^*)^{-1} \\ & = & I - (1 - |a|^2) (I - T A^*)^{-1} (I - T T^* ) (I
- A T^*)^{-1}. \een \qed

The above lemma shows in particular that the commuting tuple $T_a$
is contractive again.

\begin{Corollary}

If $T$ is a commuting contractive $n$-tuple of operators on $\clh$,
then so is $T_a$ for any $a \in \mathbb{B}_n$.

\end{Corollary}

\begin{Corollary}\label{unitary} There exists a unitary operator $U$ from
$\cld_{T_a ^*}$ to $\cld_{T^*}$ such that $U D_{T_a ^*} = D_{T^*}
S^*$ with  $S = s_a (I - TA ^*)^{-1} \in \clb(\clh)$.
\end{Corollary}

\NI \textsf{Proof.} The equality $$\| D_{T_a ^*} h \|^2 = \| D_{T
^*} S^* h \|^2$$ holds for each $h$ in $\clh$ and enables us to
define an isometry $U : \cld_{T_a ^*} \raro \cld_{T ^*}$ by $U
(D_{T_a ^*} h ) = D_{T^*} S^* h$. Since $S$ is invertible, this
isometry is even a unitary operator.  \qed

\vspace{0.3cm} Our next aim is to show that characteristic functions
behave naturally with respect to biholomorhic mappings of the unit
ball.

\begin{Definition} Let $\{ \cle , \cle_* , \varphi\}$ and $\{ \clf
, \clf_* , \psi \}$ be two operator valued bounded functions. We say
that these two functions coincide weakly if there exists a unitary
operator $\tau : \cle_* \raro \clf_*$ such that $ \psi (w) \psi
(z)^* = \tau  \varphi (w) \varphi (z)^* \tau ^*$ for all $z,w \in
\mathbb{B}_n.$
\end{Definition}

Obviously coincidence implies weak coincidence. Conversely, the
bounded operator valued functions $\{\clh,\clh,\varphi \equiv
I_{\clh}\}$ and $\{\clh \oplus \clh, \clh , \psi \equiv
(0,I_{\clh})\}$ coincide weakly, but an elementary argument shows
that they do not coincide. On the other hand, weak coincidence
almost implies coincidence.

Recall that the support of a  bounded operator valued function $\{
\cle , \cle_* , \varphi\}$, is defined as
$$
{\rm{supp}}(\varphi) = \overline{{\rm span}} \cup \{
{\rm{Ran}}\varphi(z)^* : z \in \mathbb{B}_n \} = \cle \ominus \cap
\{  {\rm{Ker}} \varphi(z) : z \in \mathbb{B}_n\}.
$$

\vspace{0.4cm}

\begin{Lemma} Let $\{ \cle , \cle_* , \varphi\}$ and
$\{ \clf , \clf_* , \psi \}$ be bounded operator valued functions.
\begin{enumerate}

\item The above functions coincide weakly if and only if their
restrictions $\{\rm{supp}(\varphi), \cle_* ,
\varphi|\rm{supp}(\varphi)\}$ and $\{\rm{supp}(\psi), \clf_* ,
\psi|\rm{supp}(\psi)\}$ coincide.

\item The functions $\{ \cle , \cle_* , \varphi\}$ and $\{ \clf ,
\clf_* , \psi \}$ coincide if and only if they coincide weakly and
$\cle \ominus \rm{supp}(\varphi)$ is isometrically isomorphic to
$\clf \ominus \rm{supp}(\psi).$ \end{enumerate}
\end{Lemma}

\vspace{0.2cm}

\NI \textsf{Proof.} Suppose that there is a unitary operator $\tau :
\cle_*  \raro \clf_*$ with
$$
\psi(w)\psi(z)^* = \tau \varphi (w) \varphi(z)^* \tau^*, \quad z,w
\in \mathbb{B}_n.
$$
Then there is a unique unitary operator $U : \rm{supp}(\varphi)
\raro \rm{supp}(\psi)$ such that
$$
U(\varphi(z)^* \tau^* x) = \psi(z)^* x, \quad z \in \mathbb{B}_n , x
\in \clf_*.
$$
Obviously this operator satisfies the intertwining relations
$$
(\psi(z)| \rm{supp}(\psi)) U = \tau (\varphi(z)|\rm{supp}(\varphi)),
\quad z \in \mathbb{B}_n.
$$

Conversely suppose that there is a unitary operator $U :
\rm{supp}(\varphi) \raro \rm{supp}(\psi)$ satisfying the last
intertwining relations. Then it is elementary to check that
$\psi(w)\psi(z)^* = \tau \varphi(w) \varphi(z)^* \tau^*$ for $z,w
\in \mathbb{B}_n.$ Furthermore, if there is a unitary operator $V:
\cle \ominus \rm{supp}(\varphi) \raro \clf \ominus \rm{supp}(\psi),$
then obviously
$$
\psi(z) (U \oplus V) = \tau \varphi(z), \quad z \in \mathbb{B}_n.
$$
To complete the proof, suppose that there is a unitary operator $W:
\cle \raro \clf$ with
$$
\psi(z) W = \tau \varphi(z), \quad z \in \mathbb{B}_n.
$$
Then necessarily $W(\bigcap_{z \in \mathbb{B}_n}{\rm Ker}
\varphi(z)) \subset \bigcap_{z \in \mathbb{B}_n} {\rm Ker} \psi(z)$,
and using the same property of $W^{-1},$ we see that $W(\cle \ominus
\rm{supp}(\varphi)) = \clf \ominus \rm{supp}(\psi).$ \qed

\vspace{0.4cm}

Any biholomorphic automorphism of the unit ball is of the form $u
\circ \varphi_a$, where $u$ is a unitary operator on $ \mathbb{C}^n$
and $a \in  \mathbb{B}_n$ (see \cite{Rudin}). Let $(u_{ij})$ be the
matrix representation of $u$. We denote by $u(T)$ the commuting
tuple $(\sum u_{1j} T_j , \ldots , \sum u_{nj} T_j)$ which is easily
seen to be contractive again. The image of $T$ under the
biholomorphic automorphism $u \circ \varphi_a$, obtained by applying
the anlytic functional calculus, is $u(T_a)$.

\bt

Let $T$ be a commuting contractive tuple of bounded operators on
$\clh$ and let $ u \circ \varphi_a$ be an arbitrary biholomorphic
automorphism of $ \mathbb{B}_n$. Then the operator valued
contractive analytic functions $ \{ \cld_{T}, \cld_{T^*}, \theta_T
\circ \varphi_a \circ u^*  \}$ and $ \{ \cld_{u(T_a) }, \cld_{u(T_a)
^*}, \theta_{u(T_a)} \}$ coincide weakly.

\et

\NI \textsf{Proof.}  It is elementary to check that the two
functions $ \{ \cld_{T}, \cld_{T^*},\theta_T \circ u^*  \}$ and $ \{
\cld_{u(T) }, \cld_{u(T) ^*}, \theta_{u(T)} \}$ coincide. Hence we
only need to prove that there is a unitary operator $U :
\mathcal{D}_{{T_a}^*} \raro \mathcal{D}_{T^*}$ such that
$$ \theta_T (\varphi_a(w)) \theta_T (\varphi_a(z))^* = U
\theta_{T_a} (w) \theta_{T_a} (z)^* U ^*$$ for $z, w \in
\mathbb{B}_n$. For $z$ in $\mathbb{B}_n$, let us abbreviate
 $\varphi_a (z)$ by $z^{\prime}$. Recall that for $z, w \in
 \mathbb{B}_n$, the identity
\begin{equation} \label{theta}
I - \theta_T(w) \theta_T(z)^* = (1 - \langle w , z \rangle) D_{T^*}
(I - W T^*)^{-1} (I - T Z^*)^{-1} D_{T^*}
\end{equation}
holds \cite[Lemma 2.2]{BES} Using the definition of $\varphi_a$ and
the observation that $(P_a w)T^* = w P_a(T)^*$, we find that \ben (I
- w^{\prime} T^*) & = & I
- \frac{a - P_a w - s_a Q_a w}{1 - \langle w , a \rangle} T^* \\
& = & (1 - \langle w , a \rangle)^{-1} ( I - w A^*  - A T^* +
(P_a w) T^* + s_a (Q_a w) T^* ) \\
& = & (1 - \langle w , a \rangle)^{-1} [ (I - A T^*) - w ( A^* -
(P_a T)^* - s_a (Q_a T)^*) ] \\& = & (1 - \langle w , a
\rangle)^{-1} [ (I - A T^*) - w T_a ^* (I - A T^*) ] \\ & = & (1 -
\langle w , a \rangle)^{-1} (I - w T_a ^*) (I - A T^*). \een By
passing to inverses we obtain that
\begin{equation} \label{iden1}(1 - \langle w , a \rangle)^{-1} (I -
w^{\prime} T^*)^{-1} = (I - A T^*)^{-1} (I - w T_a ^*)^{-1}.
\end{equation}

Replacing $w$ by $z$ leads to
\begin{equation} \label{iden2}
(1 - \langle a, z \rangle)^{-1} (I - T {z^{\prime}}^*)^{-1} = (I -
T_a z^*)^{-1} (I - T A^*)^{-1}.\end{equation}

Using equation (\ref{theta}), we see that \ben \theta_T (w^{\prime})
\theta_T (z^{\prime})^* & = & D_{T^*} - (1 - w^{\prime}
{z^{\prime}}^*) D_{T^*} (I - w^{\prime} T^*)^{-1} (I - T
{z^{\prime}}^*)^{-1} D_{T^*} \\ & = & I - (1 - \langle \varphi_a
(w), \varphi_a (z) \rangle ) D_{T^*} (I - w^{\prime} T^* )^{-1} (I -
T {z^{\prime}}^*)^{-1} D_{T^*}. \een Note that $$ 1 - \langle
\varphi_a (w), \varphi_a (z) \rangle = \frac{ (1 - |a|^2) ( 1 -
\langle w, z \rangle)}{(1 - \langle w, a \rangle) (1 - \langle a, z
\rangle)}.$$ Therefore the last equality implies that
\begin{equation} \theta_T (w^{\prime}) \theta_T (z^{\prime})^* =
I - \frac{ (1 - |a|^2) ( 1 - \langle w, z \rangle)}{(1 - \langle w,
a \rangle) (1 - \langle a, z \rangle)} D_{T^*} (I - w^{\prime} T^*
)^{-1} (I - T {z^{\prime}}^*)^{-1} D_{T^*}.
\end{equation}

On the other hand, by (\ref{theta}) and an application of Corollary
\ref{unitary}, we have \ben & & U \theta_{T_a} (w) \theta_{T_a} (z)
^* U ^* \\
& = & I - (1 - \langle w , z \rangle) U D_{T_a ^*} (I - W T_a
^*)^{-1} (I - T_a Z^*)^{-1} D_{T_a ^*} U ^* \\
& = & I - (1 - \langle w , z \rangle) D_{T^*} S^* (I - W T_a
^*)^{-1} (I - T_a Z^*)^{-1} S D_{T ^*} \\
& = & I - (1 - \langle w , z \rangle) (1 - |a|^2) D_{T^*} (I - A
T^*)^{-1} (I - W T_a ^*)^{-1} (I - T_a Z^*)^{-1} (I - T A^*)^{-1}
D_{T^*}. \een

The last equality along with (\ref{iden1})and (\ref{iden2})
completes the proof.  \qed \vspace{0.3cm}

%A few words about the module structure induced by a commuting
%contractive tuple $T$ here will be helpful for later parts of this
%note. Any such tuple $T$ makes $\clh$ a module over the polynomial
%ring $ \mathbb{C}[z_1, z_2, \ldots ,z_n]$ by means of the action
%$$ p \cdot h = p(T)h \mbox{ for } p \in \mathbb{C}[z_1, z_2, \ldots
%,z_n] \mbox{ and } h \in \clh.$$ If a bounded operator $V: \clh
%\raro \clk$ intertwines $T_i$ with $R_i$ for all $i=1,2, \ldots
%,n$, then $V$ is a homomorphism from the module $(\clh, T)$ to
%the module $(\clk , R)$. If moreover, $V$ is a unitary operator,
%then the two modules are isomorphic. A multiplier $M_\varphi$ for
%$\varphi \in M(\cle , \cle_*)$ is a module homomorphism from
%$H^2_n \otimes \cle$ to $H^2_n \otimes \cle_*$. If $E$ denotes
%the one-dimensional subspace of constant functions in $H^2_n$,
%then the action of $M_{\varphi}$ on $H^2_n \otimes \cle$ is
%determined by its action on $E \otimes \cle$ because $H^2_n
%\otimes \cle = \overline{{\rm span }} \{ (M_z^k \otimes I_{ \cle})
%(E \otimes \cle) : k \in \mathbb{N}^n\}$ and $M_{\varphi} (M_z^k
%\otimes I_\cle) = (M_z^k \otimes I_{\cle_*}) M_\varphi$ for any
%multi-index $k$.

\newsection{Model and Coincidence}

Recall that a commuting tuple $T \in \clb(\clh)^n$ is called a
spherical isometry if the column operator $\clh \raro \clh^n, \; x
\mapsto (T_i x)_{i=1}^n$ is an isometry. We shall say that $T$ is a
co-isometry if the column operator $T^*: \clh \raro \clh^n$ is an
isometry.

\vspace{0.2cm}

\begin{Definition} A commuting contractive tuple $T$ on $\clh$ is
called completely non-coisometric (c.n.c.) if there is no
non-trivial closed joint invariant subspace $\clm$ of $T_1 ^*,
\ldots, T_n ^*$ such that the tuple $P_{\clm} T |_{\clm} = (
P_{\clm} T_1 |_{\clm}, \ldots, P_{\clm} T_n |_{\clm})$ is a
co-isometry.
\end{Definition}

Given a commuting contractive $n$-tuple of operators $T$ on a
Hilbert space $\clh$, one can define a bounded operator $j : \clh
\raro H^2_n(\cld_{T^*})$ by \be j(h) = \sum_{\alpha \in
\mathbb{N}^n} \gamma_{\alpha} ( D_{T^*} T^{* \alpha} h) z^{\alpha}
\label{j}. \ee It is well known that $L = j^*: H^2_n(\cld_{T^*})
\raro \clh$ is the unique continuous linear map satisfying
$$
L (p \otimes \xi) = p (T) D_{T^*} \xi
$$ for all $p \in
\mathbb{C}[z_1, \ldots, z_n ]$ and $\xi \in \cld_{T^*}$
\cite[Theorem 4.5]{sub3}. The operator $L$ intertwines $M_{z_i}
\otimes I_{\cld_{T^*}}$ and $T_i$ for every $i=1, \ldots ,n,$ and is
closely related to the Poisson transform defined by Popescu in
\cite{Pop99}. The following lemma gives a characterization of c.n.c.
tuples in terms of its adjoint $j = L^*$.

\vspace{0.4cm}

\begin{Lemma}  The kernel of the operator $j$ is the largest
invariant subspace for $T_1 ^*, \ldots, T_n ^*$ such that $ P_{\clm}
T |_{\clm}$ is a co-isometry.
\end{Lemma}

\vspace{0.2cm}

\NI \textsf{Proof.} Since $j \, T_i ^* = ( M_{z_i} ^* \otimes
I_{\cld_{T^*}})j$, for all $i$, the kernel of $j$ is invariant for
$T_1 ^*, \ldots, T_n ^*$. From the defining formula (\ref{j}) it
follows that, for any $h \in$ Ker $j$, we have $D_{T^*} h = 0$. Now
$$ \| h \|^2 - \sum_{i = 1}^{n} \| T_i ^* h \|^2 = \langle (I -
\sum_{i = 1}^{n} T_i T_i ^*) h, h \rangle = \| D_{T^*} h \|^2 = 0.$$
Therefore the tuple $ P_{{\rm Ker} j} T |_{{\rm Ker} j}$ is a
co-isometry.

If $\clm$ is a closed subspace invariant under $T_i ^*$ for all $i =
1, \ldots, n$ such that $P_{\clm} T |_{\clm}$ is a co-isometry, then
for all $\alpha \in \mathbb{N}^n$ and $ h \in \clm$, we have
$$\| D_{T^*}
T^{* \alpha} h \|^2 = \langle (I - \sum_{i = 1}^{n} T_i T_i ^*) T^{*
\alpha} h , T^{* \alpha} h \rangle = \| T^{* \alpha} h \|^2 -
\sum_{i = 1}^{n} \| T_i ^* T^{* \alpha} h \|^2 = 0.
$$
Hence $\clm$ is contained in the kernel of $j$. \qed

\vspace{0.4cm}

Let $T$ be a commuting contractive tuple on $\clh$. The
characteristic function $\theta_T$ of $T$ induces a contractive
multiplier $M_{\theta_T} : H^2_n \otimes \cld_T \raro H^2_n \otimes
\cld_{T^*}$. More precisely, one can show that \cite[Proposition
1.2]{EandP}
$$
(1 - \theta_T (w) \theta_T (z)^*)/(1 - \langle w,z\rangle) =
k_T(w)^* k_T (z)
$$
for $z,w \in \mathbb{B}_n,$ where $k_T (z) = (I - Tz^*)^{-1}
D_{T^*}.$ The positive definiteness of the kernel on the left is
equivalent to the fact that $M_{\theta_T}$ is a contractive
multiplier.

It is well known \cite[Equation 1.11]{curv} that the intertwining
map $L$ acts as
$$
L(k(\cdot,z) \otimes x) = k_T (z)x, \quad z \in \mathbb{B}_n, x \in
\cld_{T^*},
$$
where $k: \mathbb{B}_n \times \mathbb{B}_n \raro \mathbb{C}, k(z,w)
= (1 - \langle w,z\rangle)^{-1},$ is the reproducing kernel of
$H^2_n.$ By Lemma 3.2 the tuple $T$ is completey non-coisometric if
and only if
$$
\clh = \overline{\rm span} \{k_T (z)x:  z \in \mathbb{B}_n, x \in
\cld_{T^*}.
$$
In the follwing we shall use that the dilation map $L$ and the
characteristic multiplier $M_{\theta_T}$ of $T$ satisfy the
relations \be LL^* + A_T = I_\clh. \label{LLstar} \ee \be L^* L +
M_{\theta_T} M_{\theta_T} ^* = I_{H^2_n \otimes \cld_{T^*}}.
\label{LstarL} \ee For a proof, see \cite{sub3} and \cite{BES}. In
the particular case, that $T = M_z \in L(H^2_n)^n$ is the standard
shift, the map $L$ is a unitary operator and hence $\theta_{M_z} =
0.$

To construct a functional model for a given completely
non-coisometric commuting contractive tuple $T$, let us denote by
$\Delta = (I_{H^2_n \otimes \mathcal{D}_T} - M_{\theta_T} ^*
M_{\theta_T})^{1/2}$ the defect operator of $M_{\theta_T}$.

\vspace{0.2cm}

\begin{Lemma} \label{r map}
Let $T$ be a c.n.c. commuting contractive tuple on $\clh$. Then
there is a unique contractive linear operator $ r : \clh \raro H^2_n
\otimes \cld_{T}$ such that:
\begin{enumerate}
\item $ r (k_T (z) x) = - \Delta ( k( \cdot , z) \otimes \theta_T
(z)^* x)$ for $z \in \mathbb{B}_n$ and $x \in \cld_{T^*}$;
\item $\| h \|^2 = \| j (h) \|^2 + \| r (h) \|^2$ for
all $h \in \clh$;
\item $ r L = - \Delta M_{\theta_T} ^*$.\end{enumerate}
\end{Lemma}

\NI\textsf{Proof.}  By the remarks preceding the lemma the
uniqueness is clear. To prove the existence, first observe that \ben
\lefteqn{\| \Delta \sum_{i = 1}^{m}
 k(\cdot , z^{(i)}) \otimes \theta_T (z^{(i)})^* x_i \|^2}
 \hspace{-1.2 cm} \\ & = &
  \sum_{i, j = 1}^{m} \langle (I -
 M_{\theta_T} ^* M_{\theta_T} ) (k(\cdot , z^{(i)}) \otimes \theta_T (
 z^{(i)}) ^* x_i, k(\cdot , z^{(j)}) \otimes \theta_T (z^{(j)}) ^* x_j
 \rangle \\
  & = &
 \sum_{i, j = 1}^{m} [ \frac{ \langle \theta_T (z^{(j)})
 \theta_T (z^{(i)}) ^* x_i , x_j \rangle }{1 - \langle z^{(j)},
 z^{(i)} \rangle }  -  \langle ( I - L^* L )^2 (k(\cdot , z^{(i)})
 \otimes x_i ) , k(\cdot , z^{(j)}) \otimes x_j \rangle ]
 \een
 for $z^{(1)}, \ldots , z^{(m)} \in \mathbb{B}_n$ and
 $x_1, \ldots , x_m  \in \cld_{T^*}.$
Using the identity
$$
\langle L^* L (k(\cdot , z^{(i)})
 \otimes x_i ) , k(\cdot , z^{(j)}) \otimes x_j \rangle
=  \frac{\langle  1 - \theta_T (z^{(j)})
  \theta_T (z^{(i)}) ^* x_i , x_j \rangle }{1 - \langle z^{(j)},
   z^{(i)} \rangle }
$$
we find that \ben \lefteqn{\| \Delta \sum_{i = 1}^{m}
 k(\cdot , z^{(i)}) \otimes \theta_T (z^{(i)})^* x_i \|^2}
 \hspace{-1.2 cm} \\ & = &
\sum_{i,j=1}^{m}
 [ \langle k_T (z^{(j)})^* k_T (z^{(i)}) x_i , x_j \rangle -
  \langle (L^* L)^2 k(\cdot , z^{(i)})
   \otimes x_i , k(\cdot , z^{(j)}) \otimes x_j \rangle \\ & = &
\| \sum_{i = 1}^{m} k_T (z^{(i)})
 x_i \|^2 - \| L^* ( \sum_{i =1}^{m} k_T (z^{(i)}) x_i ) \|^2.
\een Hence there is a unique contractive linear map $ r : \clh \raro
H^2_n \otimes \cld_{T}$ satisfying condition (1). The above
computation shows that condition (2) holds as well. The proof is
completed by the observation that
$$
r L (k(\cdot , z) \otimes x) = r (k_T (z) \otimes x) = - \Delta
(k(\cdot , z) \otimes \theta_T (z)^* x) = - \Delta M_{\theta_T} ^*
(k(\cdot , z) \otimes x)
$$
for all $z \in \mathbb{B}_n$ and $x \in \cld_{T^*}$. \qed

\vspace{0.3cm}

The observation (2) of the lemma above allows us to define an
isometry
$$
V: \clh \raro  (H^2_n \otimes \cld_{T^*}) \oplus \overline{{\rm Ran}
\Delta}, \; V h = jh \oplus rh$$ Our next aim is to show that the
range of $V$ is the orthogonal complement of the range of the
isometry $U: H^2_n \otimes \cld_{T} \raro (H^2_n \otimes \cld_{T^*})
\oplus \overline{{\rm Ran} \Delta}$ defined by $U \xi = M_{\theta_T}
\xi \oplus \Delta \xi$ for $\xi \in H^2_n \otimes \cld_T$.

\vspace{0.2cm}

\begin{Lemma}\label{isometry} Suppose $T$ is a c.n.c.
commuting contractive tuple.
Then the isometries $U$ and $V$ defined above satisfy the relation
$$
 V V^* + U
U^* = I_{ (H^2_n \otimes \cld_{T^*}) \oplus \overline{{\rm Ran}
\Delta}}.
$$
\end{Lemma}

\NI \textsf{Proof.} Note that the block operator matrix for $V V^* +
U U^* $ with respect to the decomposition $(H^2_n \otimes
 \cld_{T^*}) \oplus
\overline{{\rm Ran} \Delta}$ is
\begin{equation}
 \begin{bmatrix} L^* L + M_{\theta_T} M_{\theta_T} ^*& L^* r^*
  + M_{\theta_T} \Delta \\ r L + \Delta M_{\theta_T} ^*  & r r^*
  + \Delta ^2
   \end{bmatrix}
\end{equation}
We know that $L^* L + M_{\theta_T} M_{\theta_T} ^* = I_{H^2_n
\otimes \cld_{T^*}}$ and $ r L +  \Delta M_{\theta_T}^*= 0$. So all
that remains is to show that $r r^* + \Delta ^2$ is the orthogonal
projection onto $\overline{{\rm Ran }\Delta}$.

By definition, ${\rm Ran} \, r \subseteq \overline{{\rm Ran } \,
\Delta},$ and therefore ${\rm Ker } \, \Delta \subseteq {\rm Ker }
\, r^* = {\rm Ker } \, r r^*$.  Hence $r r^* + \Delta ^2$ is zero on
$({\rm Ran} \, \Delta)^{\perp}$. Using condition (1) in Lemma 3.3 we
find that \ben \lefteqn{\langle r^* \Delta (k(\cdot , z) \otimes x)
, k_T (w) y \rangle} \hspace{-1.88888888cm} \\ & = &
 - \langle \Delta ^2 (k(\cdot , z) \otimes x), M_{\theta_T}
^* (k(\cdot , w) \otimes y) \rangle
\\
 & = & - \langle M_{\theta_T}
(I - M_{\theta_T} ^* M_{\theta_T} ) (k(\cdot , z) \otimes x),
k(\cdot , w) \otimes y
 \rangle \\
& = & - \langle ( I - M_{\theta_T} M_{\theta_T}^*) M_{\theta_T}
(k(\cdot , z) \otimes x) , k(\cdot , w) \otimes y \rangle \\
& = & - \langle L M_{\theta_T} (k(\cdot , z) \otimes x) , k_T (w) y
\rangle. \een holds for all $z,w \in \mathbb{B}_n$ and $x \in \cld_T
, y \in \cld_{T^*}.$ Therefore we obtain the intertwining relation
 \begin{equation} \label{r eqn} r^* \Delta = - L
M_{\theta_T}.\end{equation} But then the observation that
$$
r r^* \Delta + \Delta ^3 = - r L M_{\theta_T} + \Delta ^3 = \Delta
M_{\theta_T} ^* M_{\theta_T} + \Delta ^3 = \Delta ( I - \Delta ^2) +
\Delta ^3 = \Delta
$$
suffices to complete the proof. \qed

\vspace{0.3cm}

Thus $V$ and $U$ are isometries with orthogonal ranges such that
$$
{\rm Ran} V \oplus {\rm Ran} U = (H^2_n \otimes \cld_{T^*}) \oplus
\overline{{\rm Ran} \Delta}.
$$

According to Lemma 3.4 the isometry $V$ induces a unitary operator
between $\clh$ and the space
$$
\mathbb{H}_T = ((H^2_n \otimes \cld_{T^*}) \oplus \overline{{\rm
Ran} \Delta}) \ominus \{ ( M_{\theta_T} u, \Delta u): u \in H^2_n
\otimes \cld_T \}.
$$
To prove that the characteristic function is a complete unitary
invariant we shall give a functional description of the operator
tuple $T$ with the help of the above unitary operator $V$. Define
$\mathbb{T}_i \in \clb(\mathbb{H}_T)$ for $i = 1, \ldots, n$ by
$\mathbb{T}_i = V T_i V^*|_{\mathbb{H}_T}$. Given $(u, v) \in
\mathbb{H}_T$, let $h$ be the vector in $\clh$ such that $(u, v) = V
h$. Then it follows that
\begin{equation} \label{model-formula}
\mathbb{T}^*_i (u, v) = \mathbb{T}^*_i V h = V T^*_i h = (M^*_{z_i}
\otimes I_{\cld_T}) jh.
\end{equation}
The vector $rT^*_i h$ is contained in $\overline{{\rm Ran}\Delta} =
(\rm{Ker}\Delta)\perp.$ Using (3.5) we see that
$$
 \Delta rT^*_i = - M^*_{\theta_T}L^* T^*_i h =
 - M_{\theta_T}^* (M^*_{z_i} \otimes I_{\cld_{T^*}}) L^* h
=  (M_{z_i}^* \otimes I_{\cld_{T^*}}) \Delta r h .
$$
So if $\Delta ^{-1} : {\rm Ran} \Delta \raro ({\rm Ker}
\Delta)^{\perp} = \overline{{\rm Ran} \Delta}$ denotes the inverse
of the bijective linear map $\Delta : ({\rm Ker} \Delta)^{\perp}
\raro {\rm Ran} \Delta,$ then $r T_i^* h = (\Delta^{-1} (M_{z_i}^*
\otimes I_{\cld_{T^*}}) \Delta) (r h).$ Thus in view of
(\ref{model-formula}), we have constructed the following functional
model for any given completely non-coisometric commuting
 contractive tuple $T.$

\vspace{0.3cm}

\begin{Theorem}  Let $T$ be a c.n.c. commuting contractive tuple on
a Hilbert space $\clh$, and let the Hilbert space $\mathbb{H}_T$ be
defined as above. Then $T$ is unitarily equivalent to the tuple
$\mathbb{T} \in \clb(\mathbb{H}_T)^n$ whose action is given by
$$\mathbb{T}_i ^* (u, v) = ( (M_{z_i}^* \otimes I_{\cld_{T^*}})u,
\Delta ^{-1} (M_{z_i}^* \otimes I_{\cld_T}) \Delta v)$$ for $(u, v)
\in \mathbb{H}_T$ and $i = 1, \ldots,n$.

\end{Theorem}

\vspace{0.3cm}

As an application of the functional model constructed above, we
prove that the characteristic function is a complete unitary
invariant for completely non-coisometric commuting contractive
tuples.

\vspace{0.3cm}

\begin{Theorem} Suppose that $T \in \clb(\clh)^n$ and $R \in
\clb(\clk)^n$ are c.n.c. commuting contractive tuples. Then the
characteristic functions $\theta_T$ and $\theta_R$ coincide if and
only if $T$ and $R$ are unitarily equivalent. \label{coincidence
implies u.e.}
\end{Theorem}

\NI \textsf{Proof.} Suppose that $T$ and $R$ are unitarily
equivalent, that is, there is a unitary operator $U : \clh \raro
\clk$ such that
$$
U T_i = R_i U \;\;\;\; (1 \leq i \leq n).
$$
Then it is elementary to prove that the operators $ \oplus U : \clh
^n \raro \clk ^n$ and $U : \clh \raro \clk$ induce unitary operators
$$
\tau = \oplus U : \cld_T \raro \cld_R \;\;\; {\rm and} \;\;\; \tau_*
= U : \cld_{T^*} \raro \cld_{R^*}
$$
such that $\theta_T$ and $\theta_R$ coincide via $\tau$ and
$\tau_*$.

Conversely suppose that there are unitary operators $\tau^{\prime} :
\cld_T \raro \cld_R$ and $\tau_* ^{\prime} : \cld_{T^*} \raro
\cld_{R^*}$ with $$\tau_* ^{\prime} \theta_T (z) = \theta_R (z) \tau
^{\prime} \;\;\;\;\;\; (z \in \mathbb{B}_n).$$ Then the induced
operators $\tau = I \otimes \tau^{\prime} : H^2_n \otimes \cld_T
\raro H^2_n \otimes \cld_R$ and $\tau_* = I \otimes \tau_* ^{\prime}
: H^2_n \otimes \cld_{T^*} \raro H^2_n \otimes \cld_{R^*}$ are
unitary and satisfy the relations $
 \tau_* M_{\theta_T} =
M_{\theta_R}
 \tau.
$
 It follows that $$\Delta_T = (I - M_{\theta_T} ^*
M_{\theta_T})^{1/2} = ( I - \tau^* M_{\theta_R} ^* M_{\theta_R}
\tau)^{1/2} = \tau ^* \Delta_R \tau.$$ Since the diagram

\setlength{\unitlength}{3mm}
 \begin{center}
 \begin{picture}(42,17)(0,0)
 \put(0,3){$(H^2_n \otimes \cld_{T^*}) \oplus \overline{\Delta_T (H^2_n \otimes \cld_T)}$}
 \put(25,3){($H^2_n \otimes \cld_{R ^*}) \oplus \overline{\Delta_R (H^2_n \otimes \cld_R)}$}
 \put(19.5,2.2){$ \tau_* \oplus \tau$}
 \put(1.9,7.0){$ \begin{pmatrix}
   M_{\theta_T} \\
   \Delta_T \
 \end{pmatrix}$}

  \put(33,7.0){$\begin{pmatrix}
    M_{\theta_R} \\
    \Delta_R \
  \end{pmatrix}$}
 \put(5,11){$ H^2_n \otimes \cld_{T}$} \put(29.5,11){$ H^2_n \otimes \cld_R$}
 \put(20.6,12){$ \tau$}
 \put(17.5,3.5){ \vector(1,0){6.5}} \put(17.5,11.5){ \vector(1,0){6.5}}
 \put(7.3,10.3){ \vector(0,-1){5.0}} \put(32,10.2){ \vector(0,-1){5.0}}
 \end{picture}
 \end{center}
commutes, we obtain the unitary operator $ \tau_* \oplus \tau :
\mathbb{H}_T \raro \mathbb{H}_R$ between the model spaces of $T$ and
$R$. We still have to prove that via this unitary operator the
functional models $\mathbb{T} \in \clb(\mathbb{H}_T)^n$ and
$\mathbb{R} \in \clb(\mathbb{H}_R)^n$ of $T$ and $R$ are unitarily
equivalent. Thus we have to prove the identity
$$((M_{z_i}^* \otimes I_{\cld_{R^*}}) \tau_* u, \Delta_R^{-1}
 (M_{z_i}^* \otimes I_{\cld_R}) \Delta_R \tau
v) = (\tau_* (M_{z_i}^* \otimes I_{\cld_{T^*}}) u, \tau \Delta_T
^{-1} (M_{z_i}^* \otimes I_{\cld_T}) \Delta_T v)$$ for all $(u, v)
\in \mathbb{H}_T$ and $i = 1, \ldots, n$. However, the equality of
the first components follows from the definition of $\tau_*$. To
prove the equality of the second components, denote by $\xi$ the
unique element in $({\rm Ker} \Delta_T)^{\perp} = \overline{{\rm
Ran} \Delta_T}$ with
$$\Delta_T \xi = (M_{z_i}^* \otimes I_{\cld_T}) \Delta_T v.$$ Then $\tau \Delta_T
^{-1} (M_{z_i}^* \otimes I_{\cld_T}) \Delta_T v = \tau \xi \in
\overline{{\rm Ran} \Delta_R}$ satisfies
$$ \Delta_R \tau \xi = \tau \Delta_T \xi = \tau (M_{z_i}^* \otimes
I_{\cld_{T}}) \Delta_T v = (M_{z_i}^* \otimes I_{\cld_R})\tau
\Delta_T v = ( M_{z_i}^* \otimes I_{\cld_R}) \Delta_R \tau v.$$ Thus
the second components also coincide.

\NI Since both $T \in \clb(\clh)^n$ and $R \in \clb(\clk)^n$ are
unitarily equivalent to their functional models $\mathbb{T} \in
\clb(\mathbb{H}_T)^n$ and $\mathbb{R} \in \clb(\mathbb{H}_R)^n$, we
conclude that $T$ and $R$ are unitarily equivalent. \qed

\vspace{0.4cm}

In the one-dimensional case, Theorem 3.6 holds under the
hypothesis that $T$ and $R$ are completely non-unitary contractions.
A straightforward multivariable generalization of this notion
would be to call a commuting contractive tuple $T \in \clb(\clh)^n$
completely non-unitary if there is no non-zero reducing subspace
$M \subset \clh$ for $T$ such that $T|M$ is a spherical unitary,
that is, a normal spherical isometry. For $n \geq 2$, the 
non-trivial implication of Theorem 3.6 does no longer hold under
the weaker hypothesis that $T$ and $R$ are completely non-unitary.
An elementary example is the following.\\

Let $V \in \clb(\clh)$ be a completely non-unitary co-isometry
on a complex Hilbert space $\clh \neq 0$ (e.g., the unilateral
left  shift). Then the commuting pairs $T = (V,0) \in \clb(\clh)^2$
and $R = (0,V) \in \clb(\clh)^2$ are completely non-unitary
commuting contractive tuples which are certainly not unitarily
equivalent. Since $D_T = D_V \oplus 1_{\clh}$,
$D_R =  1_{\clh} \oplus D_V$ and $ D_{T^*} = 0 = D_{R^*},$ 
the characteristic functions of $T$ and $R$ coincide.\\
%A remark is in order. If two commuting contractive tuples of
%operators are unitarily equivalent, then their characteristic
%functions coincide which in turn implies that the characteristic
%functions coincide weakly. Here is an easy construction of two
%tuples the characteristic functions of which coincide weakly, but
%the tuples are not unitarily equivalent. Compare this with Theorem
%\ref{weak to strong} which along with Theorem \ref{coincidence
%implies u.e.} implies that this could not have happened if the
%tuples in question were c.n.c. Let $\clh_1$ and $\clh_2$ be any two
%complex, separable Hilbert spaces, let $A = (A_1, A_2, \ldots ,
%A_n)$ be a commuting tuple of operators on $\clh_1$ satisfying $\sum
%A_i A_i^* = 1_{\clh_1}$ and let $B = (B_1, B_2, \ldots , B_n)$ be a
%commuting contractive tuple of operators on $\clh_2$. For $A$, one
%can take any non-degenerate representation of $C( \mathbb{B}_n)$ in
%which case $A_i$ are also normal. Consider the Hilbert space $\clh =
%\clh_1 \oplus \clh_2$ and the tuple $T = (T_1, T_2, \ldots ,T_n)$
%acting on it where $T_i = A_i \oplus B_i$. Because of the presence
%of $A$, the tuple $T$ is not c.n.c. It is easy to see that the
%characteristic functions of the tuples $T$ and $B$ coincide weakly.
%However, they are clearly not unitarily equivalent.

We now relate our functional model to the model constructed by
M\"{u}ller and Vasilescu in \cite{MV}.

\vspace{0.3cm}

\begin{Proposition} Given a commuting contractive c.n.c. tuple $T$
on $\clh$, there is a unique isometry $\varphi : \overline{\rm{Ran}
A_T} \raro \overline{{\rm Ran} \Delta}$ such that $r = \varphi A_T
^{1/2}.$
\end{Proposition}

\NI \textsf{Proof.} Since for all $h \in \clh$ the identity
$$
\| A_T ^{1/2}h \|^2 = \|h \|^2 - \| L^* h \| ^2 = \|h \|^2 - \|jh
\|^2 = \| rh \| ^2
$$
holds, there is a unique isometry $\varphi : \overline{\rm{Ran}A_T}
= \overline{\rm{Ran}A_T ^{1/2}} \raro \overline{\rm{Ran}\Delta}$
such that $\varphi A_T ^{1/2} =r.$  \qed

\vspace{0.3cm}

Note that, for all vectors $ h\in \clh$ the equality
$$
\sum_{i=1}^n \| A_T ^{1/2} T_i^* h \|^2 =
\langle (\sum_{i=1}^n T_i A_T T_i^* )h,h \rangle \\
= \langle P_T(A_T)h,h \rangle = \langle A_T h,h \rangle =
 \| A_T ^{1/2}h \|^2
$$
holds. Hence there are bounded operators $U_i :
\overline{\rm{Ran}A_T} \raro \overline{\rm{Ran}A_T}$ such that
$$
U_i (A_T ^{1/2} h) = A_T ^{1/2} T_i^*h
$$
for $i = 1, \ldots , n$ and $h \in \clh.$ The operators $U_i $
commute with each other and satisfy
$$
\sum_{i=1}^n \| U_i h \|^2 = \| h \|^2  \quad (h \in
\overline{\rm{Ran}A_T}).
$$

Let $W \in \clb(\overline{\rm{Ran}\Delta})^n$ be a spherical
isometry such that
$$
W_i h = \varphi U_i \varphi ^* h \quad (i = 1, \ldots ,n, \; h \in
\overline{\rm{Ran}\varphi}).
$$
Using the notation introduced earlier in this section, we obtain
that \ben \lefteqn{ [(M_{z_i}^* \otimes I_{\cld_{T^*}}) \oplus W_i]
Vh } \hspace{-1.2cm} \\ & = & (jT_i^* h, W_i rh) = (j T_i^* h, W_i
 \varphi A_T ^{1/2} h) \\ & = &
(jT_i^* h, \varphi U_i A_T ^{1/2} h) = (jT_i^* h, \varphi A_T ^{1/2}
T_i^* h) \\ & = & V T_i^* h = \mathbb{T}_i^* (Vh) \een for $ i = 1,
\ldots ,n$ and $h \in \clh .$ Therefore $\mathbb{T}^* \in \clb
(\mathbb{H}_T)^n$ is the restriction of the commuting tuple $ (M_z^*
\otimes I_{\cld_{T^*}}) \oplus W$ on $(H^2_n \otimes \cld_{T^*})
\oplus \overline{\rm{Ran}\Delta}$ to the invariant subspace
$\mathbb{H}_T.$ Summarizing we obtain the completely non-coisometric
case of a result of M\"{u}ller and Vasilescu \cite[Theorem 11]{MV}.

\vspace{0.2cm}

\begin{Theorem} Let $T$ be a c.n.c. commuting contractive tuple.
Then there is a spherical isometry $W$ on
$\overline{\rm{Ran}\Delta}$ such that $T^*$ is unitarily equivalent
to the restriction of the tuple
$$ (M_z^* \otimes I_{\cld_{T^*}}) \oplus W
\in \clb ((H^2_n \otimes \cld_{T^*}) \oplus
\overline{\rm{Ran}\Delta})^n
$$
to the invariant subspace $\mathbb{H}_T.$
\end{Theorem}

By a result of Athavale \cite[Proposition 2]{Ath} the spherical
isometry $W$ extends to a spherical unitary, that is, a commuting
tuple $N = (N_1, \ldots ,N_n)$ of normal operators satisfying the
identity $\sum_{i=1}^n N_i N_i^* = I.$

Our final aim in this section is to show that in the class of c.n.c.
commuting contractive tuples, it is enough to consider weak
coincidence of characteristic functions.

\vspace{0.2cm}

\begin{Theorem}
Let $T$ and $R$ be two c.n.c commuting contractive tuples of
operators acting on the Hilbert spaces $\clh$ and $\clk$
respectively. If the two analytic operator valued functions
$\{\cld_T, \cld_{T^*}, \theta_T\}$ and $\{\cld_R, \cld_{R^*},
\theta_R\}$ coincide weakly, then they coincide. \label{weak to
strong}
\end{Theorem}
\vspace{0.2cm}

\NI \textsf{Proof.} By definition of weak coincidence, there is a
unitary $\tau : \cld_{T^*} \raro \cld_{R^*}$ such that for all  $z,
w \in \mathbb{B}_n$, we have $\theta_R(w) \theta_R(z)^* = \tau
\theta_T(w) \theta_T(z)^* \tau^*$ and hence $(I_{\cld_{R^*}} -
\theta_R(w) \theta_R(z)^*) = \tau (I_{\cld_{T^*}} - \theta_T(w)
\theta_T(z)^*) \tau^*.$ Using (\ref{theta}), we get,
$$D_{R^*} (I - W R^*)^{-1} (I - R Z^*)^{-1} D_{R^*} = \tau D_{T^*}
(I - W T^*)^{-1} (I - T Z^*)^{-1} D_{T^*} \tau^*, \;\;\; \mbox{ for
all } z, w \in \mathbb{B}_n.$$ Now letting $k_T(z) = (I - T
Z^*)^{-1} D_{T^*}$ for all $z \in \mathbb{B}_n$, we have
$$k_R(w)^* k_R(z) = \tau k_T(w)^* k_T(z) \tau^*, \;\;\;\;\;\;
\mbox{ for all } z, w \in \mathbb{B}_n.$$ A standard uniqueness
result for factorization of operator valued positive definite maps
implies now that there is a unitary $$ \label{k-unitary} U :
\overline{\mbox{span}} \{ k_R(z) \xi : z \in \mathbb{B}_n, \xi \in
\cld_{R^*} \} \raro \overline{\mbox{span}} \{ k_T(z) \eta : z \in
\mathbb{B}_n, \eta \in \cld_{T^*} \} $$ which satisfies

\bea \label{k-inter} U k_R(z) \xi = k_T(z) \tau^* \xi .\eea Now note
that $k_R(z) \xi = L_R (k(\cdot, z) \otimes \xi)$ so that we get
from (\ref{k-inter}) that $$U L_R = L_T (I \otimes \tau^*).$$
Invoking the c.n.c assumption, we see that $\clh =
\overline{\mbox{Ran}} L_T = \overline{\mbox{span}} \{ k_T(z) \eta :
z \in \mathbb{B}_n, \eta \in \cld_{T^*} \}$ and $\clk =
\overline{\mbox{Ran}} L_R = \overline{\mbox{span}} \{k_R(z) \xi : z
\in \mathbb{B}_n, \xi \in \cld_{R^*}\}$ so that $U$ is a unitary
from $\clk$ to $\clh$. We shall show that $U R_i = T_i U$ for all $i
= 1, 2, \ldots, n$. It is enough to show that $U R_i L_R = T_i U
L_R$ for all $i = 1, 2, \ldots, n$. But

\ben U R_i L_R & = & U L_R (M_{z_i} \otimes I_{\cld_{R^*}}) = L_T (I
\otimes \tau^*) (M_{z_i} \otimes I_{\cld_{R^*}}) \\ & = &
L_T(M_{z_i} \otimes I_{\cld_{R^*}}) (I \otimes \tau^*) = T_i L_T (I
\otimes \tau^*) \\ & = & T_i U L_R.\een Hence the proof is complete.
\qed

\vspace{1cm}

\newsection{A Beurling-Lax-Halmos theorem and
characteristic functions}

\vspace{0.4cm}

A function $\varphi \in M(\cle , \cle_*)$ is called purely
contractive if $\| \varphi(0) \eta \| < \| \eta \|$ for all non-zero
$\eta \in \cle$, and it is called inner if $M_\varphi$ is a partial
isometry. The characteristic function is always purely contractive.
It is inner when the tuple is pure. The last assertion follows from
(\ref{LstarL}). Our first aim in this section is to prove the
following version of the classical Beurling-Lax-Halmos theorem (cf.
\cite{McTr}).

\bt Let $\cle$ be a Hilbert space. Then a closed subspace $\clm$ of
$H^2_n \otimes \cle$ is invariant under $M_z \otimes I_\cle$ if and
only if it is of the form
$$\clm = (H^2_n \otimes \clx) \oplus \cly,$$
where $\clx$ is a closed subspace of $\cle$ and $\cly$ is a closed
subspace of $H^2_n \otimes \cle$ which is invariant under $M_z
\otimes I_\cle$ and contains no reducing subspace of $M_z \otimes
I_\cle$. Moreover, there is a Hilbert space $\clf$ and a purely
contractive, inner function $\varphi \in \clm(\clf, \cle)$ such that
$\cly =  M_\varphi (H^2_n \otimes \clf)$. \label{invariant} \et

\vspace{0.4cm}

To prove the Theorem \ref{invariant}, we need some preparations.

\vspace{0.4cm}

\begin{Lemma} A closed subspace $\clm$ of $H^2_n \otimes \cle$ is
reducing for the multiplication tuple $M_z \otimes I_{\cle}$ if and
only if there exists a closed subspace $\cll$ of $\cle$ such that
$\clm = H^2_n \otimes \cll$. \label{reducing}
\end{Lemma}

\NI \textsf{Proof.} Obviously every subspace of the form $H^2_n
\otimes \cll$, where $\cll$ is a closed subspace of $\cle$, is
reducing for $M_z \otimes I_{\cle}$.

Conversely, let $\clm$ be a reducing subspace. Denote by $P_{\clm}$
the orthogonal projection onto $\clm$ and by $P_{\cle} \in
\clb(H^2_n \otimes \cle)$ the orthogonal projection onto the
subspace of all constant $\cle -$valued functions. Then
$$
P_{\cle} = (I_{H^2_n} - \sum_{i=1}^{n} M_{z_i} M_{z_i}^*) \otimes
I_{\cle}
$$
and hence $P_{\clm} P_{\cle} = P_{\cle} P_{\clm}.$ Define $\cll =
P_{\clm} \cle.$ Then $\cll = \clm \cap \cle \subset \cle$ is a
closed subspace with
$$
H^2_n \otimes \cll = \overline{\rm span} \{z^k \otimes h : k \in
\mathbb{N}^n \; {\rm and} \; h \in \cll \} \subset \clm.
$$
To show the opposite inclusion, let $f = \sum_{k \in \mathbb{N}^n}
z^k \otimes \eta_{k} \in \clm$ with $\eta_{k} \in \cle$ be given.
Then the proof is completed by the observation that
$$
 f = P_\clm f = \sum_{k
\in \mathbb{N}^n} P_{\clm} (z^k \otimes \eta_k ) = \sum_{k \in
\mathbb{N}^n} P_{\clm} (M_z ^k \otimes I_{\cle}) (1 \otimes \eta_k)
= \sum_{k \in \mathbb{N}^n} (M_z ^k \otimes I_{\cle}) P_{\clm} (1
\otimes \eta_k) \in H^2_n \otimes \cll.
$$
 \qed

\vspace{0.1in}

\begin{Lemma} Let $\cln$ be an invariant subspace for the
tuple $M_z \otimes I_\cle$ on $H^2_n \otimes \cle$. Then there is a
Hilbert space $\clf$ and a purely contractive inner function
$\varphi \in \clm(\clf, \cle)$ such that $\cln = M_\varphi (H^2_n
\otimes \clf)$ if and only if $\cln$ does not contain any non
trivial reducing subspace of $M_z \otimes I_\cle$.
\label{noreducing}
\end{Lemma}

\NI \textsf{Proof.} First suppose that $\cln$ does not contain any
non trivial reducing subspace. Define $T$ to be the compression of
$M_z \otimes I_\cle$ to the subspace $\cln^{\perp}$, that is,
$$T_i = P_{\cln^{\perp}} (M_{z_i} \otimes
I_{\cle})|_{\cln^{\perp}} \mbox{ for } i = 1, 2, \ldots, n.$$ Then
$T$ is a pure commuting contractive tuple.

\vspace{3mm}

Since the $C^*-$subalgebra of $\clb(H^2_n)$ generated by $M_z$ is of
the form \be C^*(M_z) = \overline{\rm span} \{M^k_z M_z^{*j} : \;
k,j \in \mathbb{N}^n \} \label{computation} \ee
(\cite[Theorem5.7]{sub3}), the space
$$
\clm = \overline{\rm span} \bigcup \{(M_z^k \otimes I_{\cle})
\cln^{\perp} :  k \in \mathbb{N}^n \}
$$
is a reducing subspace for $M_z \otimes I_{\cle}$ which contains
$\cln^{\perp}$. Therefore $\cln$ contains the reducing subspace
$\clm^{\perp}$. Thus by hypothesis $\clm^{\perp}$ is $\{ 0 \}$ and
hence $\clm = H^2_n \otimes \cle$.

On the other hand, it is elementary to check that
$$
H^2_n \otimes \cld_{T^*} =  \overline{\rm span} \bigcup \{(M_z^k
\otimes I_{\cld_{T^*}}) j \cln^{\perp} :  k \in \mathbb{N}^n \},
$$
where $j: \cln^{\perp} \raro H^2_n \otimes \cld_{T^*}$ is the
isometry associated with the pure commuting contractive tuple $T$
according to formula (\ref{j}).

Using (4.1) one can easily show that there is a unique unitary
operator $U: H^2_n \otimes \cld_{T^*} \raro H^2_n \otimes \cle$ such
that
$$
U (M^k_z \otimes I_{\cld_{T^*}}) (jh) = ( M^k_z \otimes I_{\cle})h
$$
for all $k \in \mathbb{N}^n$ and $h \in \cln^{\perp}.$ In
particular, $U ({\rm Ran} \; j) = \cln^{\perp}.$ But then $U$ is a
unitary operator that intertwines $M_z \otimes I_{\cld_{T^*}}$ and
$M_z \otimes I_{\cle}.$ By a well known commutant lifting theorem
\cite[Therorem 5.1]{BTV}, there is a multiplier $u \in \clm
(\cld_{T^*}, \cle)$ with $U = M_u$. A standard argument shows that
$u$ has to be of the form $u \equiv \tau$ for some unitary operator
$ \tau : \cld_{T^*} \raro \cle.$ Then $\varphi (z) = \tau
\theta_T(z)$ defines a purely contractive inner multiplier $\varphi
\in \clm (\cld_T,\cle)$ with
$$
\cln = [(I \otimes \tau) {\rm Ran} \; j]^{\perp} = (I \otimes \tau)
({\rm Ran} \; j)^ {\perp} = {\rm Ran} \; M_{\varphi}.
$$

Conversely, let $\varphi \in \clm (\clf,\cle )$ be a purely
contractive inner multiplier, and let $\cll \subset \cle$ be a
closed subspace such that $H^2_n \otimes \cll \subset \rm{Ran}
M_{\varphi} .$ Then, for $\eta \in \cll$, we obtain that
$$
1 \otimes \eta = P_{\cle} (1 \otimes \eta) = P_{\cle} M_{\varphi}
M^*_{\varphi} (1 \otimes \eta ) = \varphi (0) \varphi (0)^* \eta.
$$
Since $\varphi$ is purely contractive, it follows that $\cll = \{ 0
\}.$  \qed

\vspace{4mm}

As a particular case of the above lemma we obtain the following
result for characteristic multpliers.

\vspace{4mm}

\begin{Corollary} Let $T$ be a pure commuting contractive tuple of
operators on some Hilbert space $\clh$. Then ${\rm Ran} M_{\theta_T}
$ contains no non-trivial reducing subspace for the multiplication
tuple $M_z \otimes I_{\cld_{T^*}}$.
\end{Corollary}

\vspace{2mm}

\NI \textsf{Proof of Theorem \ref{invariant}.} Let $\clm$ be a
closed subspace of $H^2_n \otimes \cle$ invariant for the tuple $M_z
\otimes I_{\cle}$. By Lemma \ref{reducing}, any reducing subspace of
$H^2_n \otimes \cle$ for the multiplication tuple $M_z \otimes
I_{\cle}$ is of the form $H^2_n \otimes \cll$ for some closed
subspace $\cll$ of $\cle$. Define
$$
\clc (\clm) = \{ \cll : \cll \mbox{ is a closed subspace of } \cle
\mbox{ and } H^2_n \otimes \cll \subseteq \clm \}
$$
and $$ \clx = \overline{\rm span} \cup \{ \cll : \cll \in \clc
(\clm) \}, \;\; \cly = \clm \ominus (H^2_n \otimes \clx).$$ Clearly
$\cly$ is an invariant subspace for $M_z \otimes I_{\cle}$ which
does not contain any non-zero reducing subspace. To complete the
proof, it suffices to apply Lemma 4.3. \qed

\vspace{0.3cm}

It was shown in \cite{BES} that the characteristic function of a
pure commuting contractive tuple is purely contractive and inner. We
end this note with the converse.

\vspace{3mm}

\bt
 Let $\cle$ and $\cle_*$ be Hilbert spaces and let $\theta \in
\clm(\cle, \cle_*)$ be purely contractive and inner. Then there is a
Hilbert space $\clh$ and a pure commuting contractive tuple of
operators $T$ on $\clh$ such that the function $\{ \cle , \cle_* ,
\theta \}$ coincides weakly with $\{ \mathcal{D}_T,
\mathcal{D}_{T^*} , \theta_T\}$. Furthermore, the tuple $T \in
\clb(\clh)^n$ is uniquely determined up to unitary equivalence.

\et

\NI \textsf{Proof.} Define $\cln = \rm{Ran} M_{\theta}$ and $ T =
P_{\cln^{\perp}} M_z|\cln^{\perp}.$ As shown in the proof of Lemma
4.3, there is a unitary operator $U = I \otimes \tau: H^2_n \otimes
\cld_{T^*} \raro H^2_n \otimes \cle_*$ such that
$$
U \circ j : \cln^{\perp} \raro H^2_n \otimes \cle_*
$$
is the inclusion mapping. Here we have used the notations
established in Lemma 4.3. Using Equation (3.2) we fint that
$$
U M_{\theta_T} M_{\theta_T}^* U^* = I - Ujj^* U^* = P_{\cln} =
M_{\theta} M^*_{\theta}.
$$
By applying both sides to vectors of the form $k(\cdot, w) \otimes
x$ and by forming the scalar product with $k(\cdot,z) \otimes y$,
one obtains that
$$
\tau \theta_T (w) \theta_T (z)^* \tau^* = \theta (w) \theta (z)^*
$$
for all $z,w \in \mathbb{B}_n.$

Suppose that $R \in \clb(\clk)^n$ is a pure commuting contractive
tuple such that $\{\cle , \cle_* , \theta \}$ and $\{\cld_R ,
\cld_{R^*} , \theta_R \}$ coincide weakly. By definition there is a
unitary operator $\sigma : \cld_{R^*} \raro \cle_*$ such that
$$
\sigma \theta_R (w) \theta_R(z)^* \sigma^* = \theta (w) \theta (z)^*
, \quad z,w \in \mathbb{B}_n.
$$
By reversing the arguments from the previous paragraph we find that
$$
V M_{\theta_R} M_{\theta_R}^* V^*= M_{\theta} M_{\theta}^*,
$$ where $V = I_{H^2_n} \otimes \sigma.$ Hence $V$ induces a
unitary operator
$$
V: (\rm{Ran} M_{\theta_R})^{\perp} \raro (\rm{Ran}
M_{\theta})^{\perp}
$$
intertwining the compressions of $M_z \otimes I_{\cld_{R^*}}$ and
$M_z \otimes I_{\cle_*}$ on both spaces. Hence $R$ and $T$ are
unitarily equivalent. \qed

\vspace{0.6cm}

The above theorem shows that up to weak coincidence each purely
contractive inner function $\theta \in \clm (\cle , \cle_*)$ is the
characteristic function of a uniquely determined pure commuting
contractive tuple $T$. It would be desirable to decide when $\{ \cle
, \cle_* , \theta \}$ and $\{\cld_T , \cld_{T^*}, \theta_T \}$ even
strongly coincide. Lemma 2.5 gives at least a first answer to this
question.

\vspace{1.0cm}

{Acknowledgement:} The third named author's research work is
supported by a UGC fellowship.

\end{document}